\documentclass[12pt]{article}
\usepackage{amsxtra,amssymb,amsthm,amsmath,latexsym}

\theoremstyle{plain}

\def\R{{\mathbb R}}

\def\oH{{\overset{\circ}{H}}}
\def\oH1{{\overset{\circ}{H}\kern-.02in{}^1}}

\def\bee{\begin{equation*}}
\def\eee{\end{equation*}}
\def\be{\begin{equation}}
\def\ee{\end{equation}}

\begin{document}

\title{Global existence of solutions to nonlinear Volterra integral equations.
}

\author{Alexander G. Ramm\\
 Department  of Mathematics, Kansas State University, \\
 Manhattan, KS 66506, USA\\
ramm@math.ksu.edu\\
http://www.math.ksu.edu/\,$\sim$\,ramm
}

\date{}
\maketitle\thispagestyle{empty}

\begin{abstract}
A new method is given for proving the global existence of the solution to nonlinear Volterra integral equations. A bound on the solution is derived. The results are based on a nonlinear inequality proved 
by the author earlier.
\end{abstract}

\section{Introduction}\label{S:1}
Consider the equation:
\be\label{e1}
u(t)=f(t)+\int_0^t a(t,s,u(s))ds, \quad t\ge 0.
\ee
The problem is: 

{\em Under what assumptions on $f$ and $a(t,s,u)$ equation \eqref{e1} has a solution which
is defined on $\R_+:=[0, \infty)$?}

Many results on the theory of integral equations and many references one can find in \cite{Z}. 
Let us formulate the author's result basic for our study (see \cite{R612}, p. 105).

Let 
\be\label{e2}
g'(t)\le -\gamma(t)g(t) +\alpha(t,g(t)) +\beta(t), \quad 
t\ge 0, \quad g\ge 0; \quad g'=\frac {dg}{dt},
\ee
where $\gamma $, $\beta$ and $\alpha$ are continuous functions of $t\in \R_+$, $\alpha(t,g)\ge 0$ is a continuous non-decreasing function of $g$ on $\R_+:=[t_0,\infty)$, $t_0\ge 0$.

{\bf Lemma 1.} {\em Assume that there exists a function $\mu=\mu(t)>0$, $\mu \in C^1([t_0, \infty))$ such that
\be\label{e3}
\alpha (t, \mu^{-1}(t))+\beta(t)\le \mu^{-1}(t)\left(\gamma (t)-\mu'(t)\mu^{-1}(t)\right), \quad \forall t\ge t_0,
\ee
and 
\be\label{e4} 
\mu(t_0)g(t_0)<1.
\ee
Then any solution $g\ge 0$ to inequality \eqref{e2}
exists on $\R_+$ and
\be\label{e5} 
0\le g(t)< \mu^{-1}(t), \quad \forall t\ge t_0.
\ee
If  $\mu(t_0)g(t_0)\le 1$ then $0\le g(t)\le \mu^{-1}(t)$ for all $t\ge t_0$.}

A proof of Lemma 1 is given in \cite{R612}, pp. 105-107, see also \cite{R657}.

 A new idea
in this paper is to use Lemma 1 with $\gamma(t)=0$.
In this case inequality \eqref{e3} may hold {\em only if $\mu(t)$ decays as $t$ grows}, and
estimate \eqref{e5} becomes the estimate of the 
{\em rate of growth of $u$}. 

In  \cite{R657} $\mu(t)$
was growing to infinity as $t\to \infty$ and estimate \eqref{e5} gave results on {\em the stability and large-time
behavior of $g(t)=\|u(t)\|$}, where the norm was a Hilbert space norm.

Let us assume that $t_0=0$ and
\be\label{e6}
|f(t)|+|f'(t)|\le c_0e^{-b_0t},\quad \forall t\ge 0,
\ee
\be\label{e7}
|a(t,t,u)|\le c_1e^{-b_1t}(1+|u|^{2p}), \quad p>0,
\ee
\be\label{e8}
\int_0^t |a_t(t,s,u(s))|ds\le c_2e^{-bt}(1+|u(t)|^{2p}), \quad a_u(t,s,u)\ge 0, \quad a_t=\frac{\partial a}{\partial t}.
\ee
Assume also that $|a|+a_u\le c(R_1, R_2)$ for $t\le R_1, s\le R_1$ and $|u|\le R_2$, $a$ and $a_t$ are smooth functions of their arguments. 
 
 This assumption allows one
to use the contraction mapping principle if $t>0$ is  sufficiently small and establish the existence and uniqueness of the local solution to equation \eqref{e1}.

Differentiate \eqref{e1} with respect to $t$ and get
\be\label{e9}
u'=f'+a(t,t, u(t))+\int_0^t a_t(t,s,u(s))ds.
\ee
{\bf Lemma 2.} {\em Let $u(t)\in H$, where $H$ is a Hilbert space, $\|u\|^2=(u,u)$, $\|u\|'=\frac {d\|u\|}{dt}$. If $u(t)\in C^1(\R_+; H)$ then
\be\label{e10}
\|u\|'\le \|u'\|.
\ee
 If $u(t)\in C^1(\R_+)$, then
\be\label{e10a}
|u|'\le |u'|.
\ee
}
{\em Proof of Lemma 2.} One has $\|u\|^2=(u,u)$. Thus,
$2\|u\|' \|u\|=(u',u)+(u,u')\le 2\|u'\|\|u\|$.
Since $\|u\|\ge 0$, one gets \eqref{e10}.

If $u(t)\in C^1(\R_+)$, then $|u(t+h)|-|u(t)|\le
|u(t+h)-u(t)|$.  Divide this inequality by $h>0$ and let $h\to 0$. This yields \eqref{e10a}.
 \hfill$\Box$

Taking the absolute value of \eqref{e9}, using \eqref{e10} and setting $g(t)=|u(t)|$, one obtains
\be\label{e11}
g'\le c_0e^{-b_0t}+c_1e^{-b_1t}(1+g^{2p}(t))+c_2e^{-bt}(1+g^{2p}(t)).
\ee 

{\bf Theorem 1.} {\em If \eqref{e4} and \eqref{e6}--\eqref{e8} hold, then the solution to \eqref{e1}
exists on $\R_+$, is unique and
\be\label{e12}
|u(t)|\le ce^{qt}, \quad q>0,
\ee
 where $ q>0$ is a fixed number and $c>0$ is a sufficiently large constant. }

In Section 2 a proof of Theorem 1 is given. From this proof one can get an estimate for the constant $c$.
\section{Proof of Theorem 1}\label{S:2}

Let us apply to \eqref{e11} Lemma 1. Choose 
\be\label{e13} 
\mu=c_3e^{-qt}, \quad q=const >0.
\ee
Since $\mu(0)=c_3$ inequality \eqref{e4} holds if 
\be\label{e14}
g(0)c_3<1.
\ee
This inequality holds if $c_3$ is sufficiently small. 

Let $t_0=0$. One has $g(0)=|f(0)|\le c_0$. 
Thus, \eqref{e14} holds if $c_3<c_0^{-1}$. 

Inequality \eqref{e3} holds if
\be\label{e15}
\begin{split}
c_0e^{-b_0t}+c_1e^{-b_1t}+c_2e^{-bt}+c_1e^{-b_1t}(c_3e^{-qt})^{-2p}
+c_2e^{-bt}(c_3e^{-qt})^{-2p}\le\\ qc_3e^{-qt} (c_3e^{-qt})^{-2}=qc_3^{-1}e^{qt}.
\end{split}
\ee
In our case $\gamma(t)=0$,  $\beta(t)$ is the sum of the three first terms in \eqref{e15} and $\alpha$
is the sum of the fourth and fifth terms in  \eqref{e15}.

Inequality \eqref{e15} holds if $c_3$ is sufficiently small and $q>0$. 

From Lemma 1 the global existence of the solution $g\ge 0$ to \eqref{e11} follows  and  estimate \eqref{e5} yields
\be\label{e16}
|u(t)|\le e^{qt}/c_3\quad \forall t\ge 0.
\ee
The local existence and uniqueness of the solution to \eqref{e1} follows from the contraction mapping principle.
 Theorem 1 is proved. \hfill$\Box$

 {\bf Remark 1.} Without some assumptions on $f$ and $a(t,s,u)$ the solution to \eqref{e1} may not  exist globally.
 
{\em Example 1.} Let $u=1+\int_0^t u^2(s)ds$. Then $u'=u^2,$ $u(0)=1$.
A simple integration yields $u=(1-t)^{-1}$. So, the solution tends to
infinity as $t\to 1$.

{\bf Remark 2.} The method developed in this paper can be used for other decay assumptions, for example, power decay of $f$ and $a(t,s,u)$ as $t\to \infty$.  

 One may look for the $\mu(t)=c_4 (1+t)^{-r}, r>0$,
where $c_4=const$. If $c_4$ is sufficiently small then
inequality  \eqref{e5} yields
$|u(t)|\le (1+t)^{r}/c_4$, $\forall t>0$.

.


\end{document}